\documentclass[11pt,twoside,reqno]{amsart}

\usepackage{amsmath}
\usepackage{amsthm}
\usepackage{amsfonts,amssymb}
\usepackage{bbm}
\usepackage{latexsym}
\usepackage{mathrsfs}
\usepackage[babel=false]{csquotes}
\usepackage[all]{xy}
\usepackage{url}
\usepackage[pdftex]{graphicx}
\usepackage{float}
\usepackage{tikz-cd}
\usepackage{hyperref}
\usepackage{cleveref}

\newtheoremstyle{TheoremNum}
		{\topsep}{\topsep}              %%% space between body and thm
		{\itshape}                      %%% Thm body font
		{}                              %%% Indent amount (empty = no indent)
		{\bfseries}                     %%% Thm head font
		{}                             %%% Punctuation after thm head
		{ }                             %%% Space after thm head
		{\thmname{#1}\thmnote{ \bfseries #3}}%%% Thm head spec

\theoremstyle{TheoremNum}
\newtheorem{teoremaRepetido}{Theorem}

\theoremstyle{plain}
\newtheorem{lema}{Lemma}[section]

\newtheorem{teorema}[lema]{Theorem}
\newtheorem{conjetura}[lema]{Conjecture}

\theoremstyle{remark}

\theoremstyle{definition}

\crefname{teorema}{Theorem}{Theorems}
\crefname{lema}{Lemma}{Lemmas}
\crefname{proposicion}{Proposition}{Propositions}
\crefname{corolario}{Corollary}{Corollaries}
\crefname{definicion}{Definition}{Definitions}
\crefname{observacion}{Observation}{Observations}
\crefname{conjetura}{Conjecture}{Conjectures}
\crefname{ejemplo}{Example}{Examples}

\newcommand{\Z}{\mathbb{Z}}

\MakeOuterQuote{"}

\begin{document}

\title[On a generalization of the seating couples problem]{On a generalization of the seating couples problem}

\author{Daniel Kohen}
\address{Departamento de Matem\'atica, Facultad de Ciencias Exactas y Naturales, Universidad de Buenos Aires and IMAS, CONICET, Argentina}
\email{dkohen@dm.uba.ar}
\thanks{DK was partially supported by a CONICET doctoral fellowship}

\author{Iv\'an Sadofschi Costa}
\address{Departamento de Matem\'atica, Facultad de Ciencias Exactas y Naturales, Universidad de Buenos Aires}
\email{ivansadofschi@gmail.com}

\begin{abstract}
We prove a conjecture of Adamaszek generalizing the seating couples problem to the case of $2n$ seats. Concretely, we prove that given  a positive integer $n$ and $d_1,\ldots,d_n\in(\Z/2n)^*$ we can partition $\Z/2n$  into $n$ pairs with differences $d_1,\ldots,d_n$.
\end{abstract}



\maketitle

\maketitle

\section{Introduction}

Preissmann and Mischler \cite{PreissmannMischler} proved the following, confirming a conjecture of R. Bacher.

\begin{teorema}\label{TeoRey}
Let $p=2n+1$ be an odd prime. Suppose we are given $n$ elements $d_1,\ldots,d_n\in (\Z/p)^*$. Then there exists a partition of $\Z/p - \{0\}$ into pairs with differences $d_1,\ldots,d_n$.
\end{teorema}

A  simpler proof of this thoerem can be found in  \cite{KohenSadofschi}. Karasev and Petrov, independently, gave a proof of this result along  the same lines and provided further generalizations in \cite{KarasevPetrov}.
In this work, they also conjectured  two generalizations of \cref{TeoRey}, replacing $p$ by an arbitrary integer $N$. The conjecture in the case that $N$ is even is originally due to Adamaszek.

\begin{conjetura}[{\cite[Conjecture 1]{KarasevPetrov}}]\label{ConjImpar}
Let $N = 2n + 1$ be a positive integer. Suppose we are given $n$ elements $d_1,\ldots,d_n\in (\Z/N)^*$. Then there exists a partition of $\Z/N - \{0\}$ into pairs with differences $d_1,\ldots,d_n$.
\end{conjetura}

We will prove the  conjecture when $N$ is even:
\medskip

\begin{teoremaRepetido}[\ref{TeoPar}]{\normalfont ({\cite[Conjecture 2]{KarasevPetrov}}){\bfseries .}} 
Let $N=2n$ be a positive integer. Suppose we are given $n$ elements $d_1,d_2,\ldots,d_n \in (\Z/N)^*$. Then there exists a partition of $(\Z/N)$ into pairs with differences $d_1, d_2,\ldots,d_n$.
\end{teoremaRepetido}

\medskip

While finishing this paper we found out that, in his master's thesis \cite{Mezei}, T.R. Mezei suggests a possible
way to solve the conjecture that is similar to ours. Furthermore, he shows that Theorem \ref{TeoPar} holds whenever $N=2p$ for $p$ a prime number.

\section{The even case}

We recall the following version of the Cauchy-Davenport theorem.
\begin{teorema}[{\cite[1.4]{Granville}}]\label{CauchyDavenport}
If $A$ and $B$ are nonempty subsets of $\Z/N$ where $0\in B$, and $gcd(b,N)=1$ for all $b\in B\setminus\{0\}$, then
$$|A+B|\geq \min\{N,|A|+|B|-1\}.$$
\end{teorema}

Suppose that we have a partition as  in \cref{TeoPar}. Since the $d_{i}$ are odd numbers, each pair contains exactly one even number. Therefore, if \cref{TeoPar} holds there exists signs $s_i$ such that $s_1d_1+\ldots+s_nd_n\equiv 1-2+3-\ldots +(2n-1)-2n\equiv n\bmod N$.

\begin{teorema}\label{TeoSuma}
Let $N=2n$ and let $d_1,\ldots,d_n\in (\Z/N)^*$. Then there exists $s_1,\ldots,s_n\in \{1,-1\}$ such that 
$$s_1d_1+\ldots+s_nd_n\equiv n\bmod {2n}$$

\begin{proof}
It is enough to prove that there exists $I\subset\{1,\ldots,n\}$ such that 
$$\sum_{i \in I} 2d_{i}  \equiv d_{1}+d_{2}+  \cdots +d_{n} + n \bmod{2n}$$
Since $d_i$ is odd for every $i$, $d_{1}+d_{2}+  \cdots+ d_{n} + n$ is even and therefore our task is equivalent to finding $I$ such that $$\sum_{i \in I} d_{i}  \equiv \frac{d_{1}+d_{2}+  \cdots +d_{n} + n}{2} \bmod{n}.$$
Let $A_{i}=\left\lbrace d_{i},0 \right\rbrace$. Applying \cref{CauchyDavenport} inductively, we see that 
$$ \#(A_{1}+  \cdots + A_{n}) \ge min \left\lbrace n, \sum \#A_{i} -(n-1) \right\rbrace=n,$$ concluding the proof.
\end{proof}
\end{teorema}

The last ingredient is the following theorem by Hall.
\begin{teorema}[{\cite{Hall}}]\label{TeoHall}
Let $A$ be an abelian group of order $n$ and $a_1,\ldots,a_n$ be a numbering of the elements of $A$. Let $d_1,\ldots,d_n\in A$ be elements such that $d_1+\ldots+d_n=0$. Then there are permutations $\sigma,\tau\in S_n$ such that
$$a_i-a_{\sigma(i)}=d_{\tau(i)}$$
\end{teorema}

\begin{teorema}\label{TeoPar}
Let $N=2n$ be a positive integer. Suppose we are given $n$ elements $d_1,d_2,\ldots,d_n \in (\Z/N)^*$. Then there exists a partition of $\Z/N$ into pairs with differences $d_1, d_2,\ldots,d_n$.
\begin{proof}
First from \cref{TeoSuma}, we may assume that $d_1+\ldots+d_n\equiv n \bmod{2n}$. Now it is enough to find a numbering $a_1,\ldots,a_n$ of the odd numbers in $\Z/N$ and $\sigma \in S_n$ such that $2i-a_i\equiv d_{\sigma(i)}\bmod {2n}$ for every $i=1,\ldots,n$, for then the partition in pairs $\{2,a_1\},\{4,a_2\},\ldots,\{2n,a_n\}$ works.

Equivalently, we need to find a numbering $b_1,\ldots,b_n$ of the even numbers in $\Z/N$ such that $2i-b_i\equiv d_{\sigma(i)}+1\bmod {N}$ for some $\sigma \in S_{n}$. Now since $d_i+1$ is even for all $i$, this is the same as finding a permutation $c_1,\ldots,c_n$ of $\{1,\ldots,n\}$ such that
$i-c_i\equiv \frac{d_{\sigma(i)}+1}{2} \bmod n$, for some $\sigma \in S_n$. 

If we verify that $\frac{d_1+1}{2}+\ldots+\frac{d_n+1}{2}\equiv 0\bmod n$ this will follow from \cref{TeoHall}. But this holds, since $d_1+\ldots+d_n\equiv n\bmod {2n}$ and therefore $(d_1+1)+\ldots+(d_n+1)\equiv 0\bmod {2n}$, proving that $\frac{d_1+1}{2}+\ldots+\frac{d_n+1}{2}\equiv 0\bmod n$.
\end{proof}
\end{teorema}

\bibliographystyle{plain} 
\bibliography{ref}

\end{document}